\documentclass[11pt]{article}
\usepackage[applemac]{inputenc}

\usepackage{amssymb}
\usepackage{amsmath}
\usepackage{graphics}
\usepackage{latexsym}
\usepackage{ae,aecompl}
\usepackage{amsthm}
\usepackage{mathrsfs}
\usepackage{amsfonts}
 \usepackage[pdftex]{graphicx}
\newcommand{\sun}{\mathbb{S}_{1}}
\newtheorem{theorem}{Theorem}[]
\newtheorem{proposition}[theorem]{Proposition}
\newtheorem{lemma}[theorem]{Lemma}
\newtheorem{corollary}[theorem]{Corollary}

\usepackage{enumerate}

\title{\Huge{Strong convergence of partial match queries in random quadtrees}}
\author{Nicolas Curien\\ \textit{\'Ecole Normale Sup\'erieure}}
\date{}
\begin{document}
\maketitle

\abstract{We prove that the rescaled costs of partial match queries in a random two-dimensional quadtree converge \emph{almost surely} towards a random limit which is identified as the terminal value of a martingale. %This sharpens the recent result of Broutin, Neininger and Sulzbach \cite{BNS} where a convergence in distribution was shown. 
Our approach shares many similarities with the theory of self-similar fragmentations.}
\section{Introduction} \label{introduction}

The quadtree structure is a storage system designed for retrieving multidimensional data. It has first been introduced  by  Finkel \& Bentley \cite{FB74} and was  studied thoroughly in computer science. The goal of this work is to study fine properties of the so-called partial match queries in random (uniform) two-dimensional quadtrees.

 Let us briefly recall  the model. Consider  a Poisson point process $\Pi$ on $\mathbb{R}_+ \times [0,1]^2$ with intensity $\textrm{d}t \otimes \textrm{d}x \textrm{d}y$. Let $((\tau_i,x_i,y_i), i\geq 1)$ be the atoms of $\Pi$ ranked in the increasing order of their $\tau$-component. We define a process $(\mathrm{Quad}(t))_{t\geq 0}$ with values in finite coverings of $[0,1]^2$ by closed rectangles with disjoint interiors as follows. We initially start with the unit square $ \mathrm{Quad}(0) := [0,1]^2$. At each time an atom $(\tau_i,x_i,y_i)$ of the Poisson process $\Pi$ falls in a rectangle of $ \mathrm{Quad}(\tau_{i}^-)$ it splits this rectangle into four subrectangles according to the horizontal and vertical coordinates of $x_i$ and $y_i$. Observe that a.s., for every $i \geq 1$, there exists a unique  rectangle of $\mathrm{Quad}(\tau_i)$ such that $(x_{i+1}, y_{i+1})$ is in its interior, hence the process $(\mathrm{Quad}(t))_{t\geq 0}$ is a.s.\,well defined. In this work, we chose to focus on the continuous time version of the random quadtree but all the results can be transferred to the random quadtree with a fixed number of points by standard depoissonization techniques, see e.g.\,\cite{BNSbis} or \cite[Lemma 1]{CJ10}.
 
 \begin{figure}[!h]
\begin{center}
\includegraphics[width=15cm]{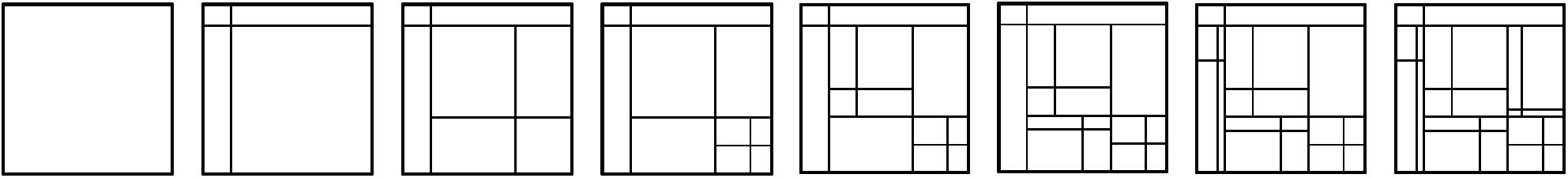}
\caption{The first $7$ splittings of a quadtree.}
\end{center}
\end{figure}

 We shall be interested in the so-called partial match query (see \cite[p 523]{FS09}).  Equivalently, for $x \in [0,1]$, we focus on the number $N_t(x)$ of rectangles in the quadtree at time $t$ whose horizontal coordinate intersects $x$, that is,
 \begin{eqnarray*}
N_t(x) & := & \# \big \{ R \in \mathrm{Quad}(t) : R \cap [(x,0),(x,1)] \neq \emptyset \big \}.
\end{eqnarray*}
The first study of the partial match has been carried out by Flajolet, Gonnet, Puech and Robson in \cite{FGPR93}. They proved that if $U$ is uniformly distributed over $[0,1]$ and independent of $ \Pi$ then  
$E[N_t(U)] $ is asymptotically equivalent to $\kappa \cdot t^{\beta}$ as $t$ tends to infinity, where \begin{eqnarray*}  \beta &:=& \displaystyle \frac{\sqrt{17}-3}{2}, \end{eqnarray*} and $\kappa >0$ is some explicit constant. The asymptotic of the expected value of $N_t(x)$ for a fixed point $x \in [0,1]$ has recently been obtained in \cite{CJ10}, it reads   \begin{eqnarray} \label{CJ} t^{-\beta}E[N_t(x)] &\xrightarrow[t\to\infty]{}& K_0\cdot h(x) ,  \end{eqnarray} where 
 \begin{eqnarray*}  K_{0}:=\displaystyle{\frac{\Gamma(2\beta+2) \Gamma(\beta+2) }{ 2 \Gamma^3(\beta+1) \Gamma^2({\beta}/{2}+1)}} \quad \mbox{ and }\quad h : u\in [0,1]  \longmapsto \big( u(1-u)\big)^{\beta/2}.  \end{eqnarray*} 
 In a very recent breakthrough \cite{BNS},  Broutin, Neininger and Sulzbach used the ``contraction method'' to obtain a convergence in distribution as $t \to \infty$ of the rescaled processes $ \{t^{-\beta}N_t(x): 0 \leq x \leq 1\}$ towards a random continuous process $\{\tilde{M}_{\infty}(x) : 0 \leq x \leq 1\}$ characterized by a recursive decomposition. The main result of the present work is to show that this convergence actually holds in a stronger sense:

\begin{theorem} \label{main} For every $x \in [0,1]$ %There exists a continuous process $(M_\infty(x))_{x \in [0,1]}$ (which is $\beta-\varepsilon$ Hölder continuous for all $\varepsilon >0$) such that 
we have the following \emph{almost sure} convergence
 \begin{eqnarray*} t^{-\beta} N_t(x)  & \xrightarrow[t\to\infty]{a.s.} & K_0 \cdot \tilde{M}_\infty(x).  \end{eqnarray*} 
\end{theorem}
\begin{figure}[!h]
 \begin{center}
 \includegraphics[width=15cm]{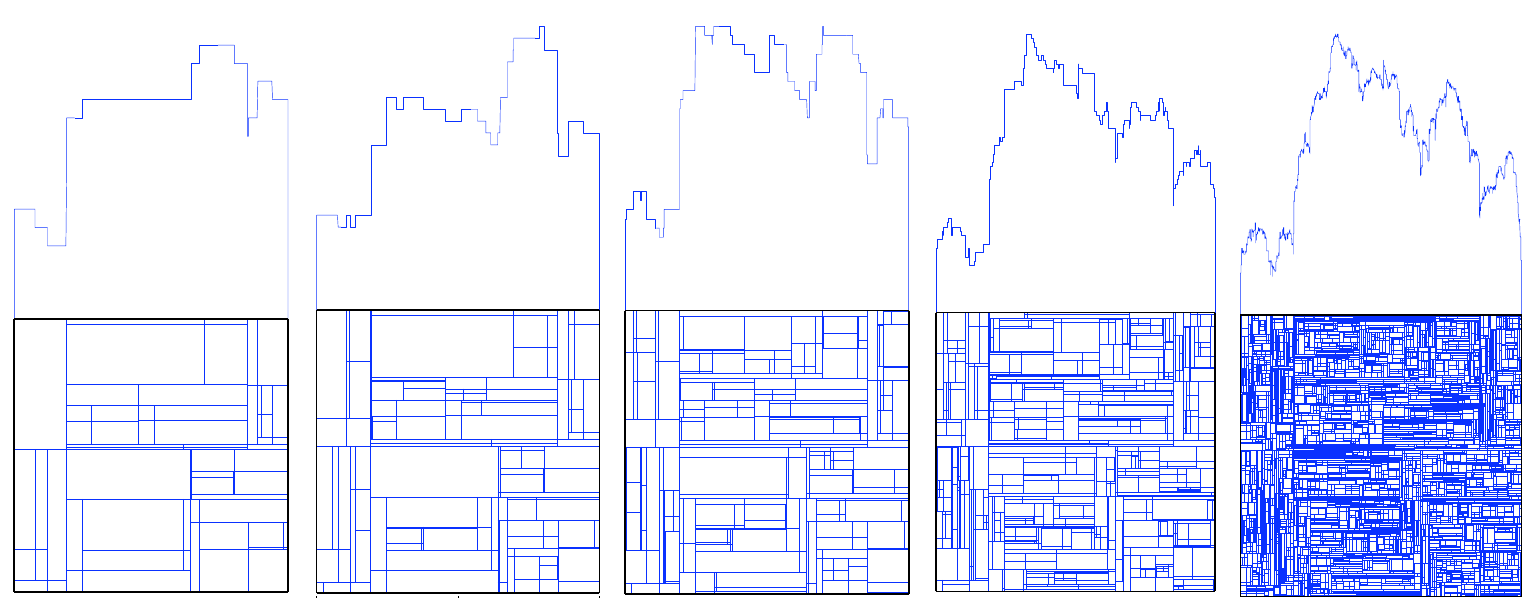}
 \caption{An illustration of the strong convergence of the partial match queries. The curves above the quadtree represent the renormalized processes $t^{-\beta} (N_{t}(x))_{x \in [0,1]}$ for $t=20,50,100,500$ and $3000$.}
 \end{center}
 \end{figure}
The random variable $\tilde{M}_{\infty}(x)$ is obtained as the limiting value of a continuous-time martingale $\{M_t(x) : t\geq 0\}$ defined by \eqref{martingale} which is a variation on the martingale introduced in \cite[Section 3.2]{BNS}. The martingale $M_t(x)$ must be considered as an analogous of the well-known Malthusian martingale in fragmentation theory, see \cite{Ber06}. Indeed, the convergence of $t^{-\beta}N_{t}(x)$ towards the limiting value of $M_{t}(x)$ is similar to the proof of \cite[Theorem 5]{BG04} and requires  some of the estimates of \cite{BNS}. See Section \ref{comments} for comments.

Here is a direct corollary of Theorem \ref{main} combined with the results of \cite{BNS}:
\begin{corollary} \label{coro}
We have the following convergence in probability \begin{eqnarray*}
\big( t^{-\beta} N_{t}(x)\big)_{x \in [0,1]}& \xrightarrow[t\to\infty]{(P)} & K_{0} \cdot \big(\tilde{M}_{\infty}(x)\big)_{x \in [0,1]},  \end{eqnarray*} for the uniform metric $\|.\|_\infty$.
\end{corollary}

The note is organized as follows: We first introduce the martingales whose limit value furnishes the process $\{ \tilde{M}_{\infty}(x) : 0 \leq x \leq 1\}$ and recall some of its properties. The third section is devoted to an estimate on the smallest and the largest rectangle in the quadtree at time $t>0$ which is used in the proof of the main result. In the last section we give some comments related to fragmentation theory.\\

\noindent \textbf{Acknowledgment: } I am grateful to Nicolas Broutin, Ralph Neininger and Henning Sulzbach for keeping me informed about their recent work on quadtrees. Special thanks go to Adrien Joseph and Henning Sulzbach for precious comments on a first version of this work.

%Throughout  this work we write 
%$$ \beta^* = \frac{\sqrt{17}-3}{2}.$$
%We aim at proving the following result: \begin{theoreme} \label{theonumeroun} For every $x \in [0,1]$,  \begin{eqnarray*} \lim_{t \rightarrow \infty} t^{-\beta^\ast} \E[N_t(x)] & = &   K_{0} \big(x(1-x)\big)^{\beta^\ast/2}, \end{eqnarray*} where the constant $K_{0}$ equals \begin{eqnarray*} K_{0}&=&\frac{\Gamma\left(2\beta^*+2\right) \Gamma(\beta^*+2) }{ 2 \Gamma^3(\beta^*+1) \Gamma^2\left(\frac{\beta^*}{2}+1\right)}. \end{eqnarray*}  \end{theoreme}

\section{The martingales} In this section we introduce the martingale  which the proof of Theorem \ref{main} is based on and compare it with the one introduced in \cite[Section 3.2]{BNS}. We start by setting some notation. \medskip

Recall the definitions of $\beta$ and of the map $h$ given in the {Introduction}. The genealogy of the rectangles appearing in the quadtree process $( \mathrm{Quad}(t))_{t\geq 0}$ can be encoded on the full infinite $4$-ary tree  \begin{eqnarray*} \mathcal{T}_4 &:=& \bigcup_{n\geq 0}\{1,2,3,4\}^n.  \end{eqnarray*} The first square $[0,1]^2$ corresponds to the word $\varnothing \in  \mathcal{T}_4$ and when a rectangle encoded by a word $u \in \mathcal{T}_4$ is split, we encode the four resulting subrectangles by $u1,u2,u3$ and $u4$ in counter clockwise order starting with the north-east rectangle. This genealogy induces a notion of ancestor, offspring... on the rectangles of $ \cup \mathrm{Quad}(t)$. The generation of a rectangle $R$ that appears in the quadtree process is the length of its encoding word in $ \mathcal{T}_4$ and is denoted by $ \mathrm{Gen}(R)$ (the length of $\varnothing$ is $0$ by convention). The time of appearance of $R$ is the first $t>0$ such that $R \in \mathrm{Quad}(t)$ and is denoted by $ \mathrm{Time}(R)$.

For $t >0$ and $x \in [0,1]$ we denote by $ \mathcal{Q}_t(x) := \{ Q_t^i(x)\}_{i \geq 1}$ the set of all rectangles belonging to $ \mathrm{Quad}(t)$ whose first coordinate intersects $x$.  The left-most and right-most horizontal coordinates of $Q_t^i(x)$ are denoted by $\ell_{t}^i$ and $r_t^i$ and we write 
  \begin{eqnarray*} x_t^i &:=& \frac{x - \ell_t ^i }{r_t^i-\ell_t^i},  \end{eqnarray*} for the ``position" of  $x$ inside the rectangle $Q_t^i(x)$. By standard properties of the Poisson point process $\Pi$, conditionally on the sigma-field $ \mathcal{F}_t$ generated by $(\mathrm{Quad}(u))_{0 \leq u \leq t}$ the number of rectangles in $ \mathcal{Q}_{t+s}(x)$ whose ancestor at time $t$ is $Q_t^i(x)$ has the same distribution as 
  $$  N'_{s\cdot \lambda(Q_t^i(x))}(x_t^i),$$ where $\lambda$ stands for the two-dimensional Lebesgue measure on $[0,1]^2$ and $N'_.(.)$ is an independent copy of the process $N_.(.)$.
 \begin{proposition} \label{propmartingale} For every $x \in [0,1]$, the process
  \begin{eqnarray} t \quad \longmapsto \quad M_t(x) &:=& \sum_{i \geq 1}  \lambda\big(Q_t^i(x)\big) ^\beta h \left( x_t^i \right),  \label{martingale} \end{eqnarray}
  is a continuous-time non-negative martingale whose limiting value is denoted by $ \tilde{M}_{\infty}(x)$. 
\end{proposition}
Before going into the proof of Proposition \ref{propmartingale} let us emphasize the difference between this martingale and the one considered in \cite{BNS}. Fix a generation $n \geq 0$ and denote by $\{\tilde{Q}_{n}^i(x)\}_{1 \leq i \leq 2^n}$ the rectangles at generation $n$ that are above the point $x \in [0,1]$ and write $\tilde{x}_{n}^i$ for the position of $x$ inside $ \tilde{Q}_{n}^i(x)$. Then from \cite[Section 3.2]{BNS}, \begin{eqnarray*}
\tilde{M}_{n}(x) &:=& \sum_{i=1}^{2^n} \lambda\big( \tilde{Q}_{n}^i(x)\big)^\beta h(\tilde{x}_{n}^i),  \end{eqnarray*}is a discrete-time non-negative martingale whose limiting value is $ \tilde{M}_{\infty}(x)$. The proof of this fact is based on the following lemma that shows that the expectation is kept after one splitting.
\begin{lemma} \label{split}We have 
 \begin{eqnarray*}E\left[ \sum_{i=1}^2 \lambda( \tilde{Q}_{1}^i(x))^\beta h(\tilde{x}_{1}^i)\right] &=& h(x). \end{eqnarray*} 
\end{lemma} 
This lemma was proved in \cite{BNS} but is also rigorously equivalent to the fact $h$ solves the integral equation that was already considered in \cite[Section 5]{CJ10}. 

The difference between the martingales $\tilde{M}_{n}(x)$ and $M_{t}(x)$ is that in the latter case we consider the splittings chronologically as they occur whereas in the first case we consider them generation after generation. It should be clear that the order in which the splittings are considered does not change the martingale property and we could use Lemma \ref{split} to show directly that $M_{t}(x)$ is a martingale for every $x \in [0,1]$. However it will be useful for our purpose to link $M_{t}(x)$ to its discrete time analog $ \tilde{M}_{n}(x)$.

\proof[Proof of Proposition \ref{propmartingale}]  Fix $x \in [0,1]$. It easily follows from \cite[Section 3.2]{BNS} (see also \cite{BNSbis}) that $\tilde{M}_{n}(x)$ converges in $L^p$ for any $p>1$ towards $\tilde{M}_{\infty}(x)$ and thus $E[\tilde{M}_{\infty}(x) \mid \mathcal{F}_{t}] = \lim E[\tilde{M}_{n}(x) \mid \mathcal{F}_{t}]$ almost surely as $n \to \infty$. By the Markov property applied at time $t>0$ and using the martingale structure of $\tilde{M}_{n}(x)$ we deduce that 
 \begin{eqnarray*} E[\tilde{M}_{n}(x) \mid \mathcal{F}_{t}] &=& \sum_{i \geq 1} \lambda(Q_{t}^i(x))^\beta h(x_{t}^i) \mathbf{1}_{ \mathrm{Gen}(Q_{t}^i(x)) < n} + \sum_{i = 1}^{2^n} \lambda(\tilde{Q}_{n}^i(x))^ \beta h(\tilde{x}_{n}^i) \mathbf{1}_{ \mathrm{Time}(\tilde{Q}_{n}^i(x)) \leq t}. \end{eqnarray*}
It is easy to see that $\inf \{\mathrm{Time}(\tilde{Q}_{n}^i(x)): 1 \leq i \leq 2^n\}$ goes to infinity as $n \to \infty$ a.s., hence letting $n$ tend to infinity in the last display we get that 
$E[\tilde{M}_{\infty}(x)\mid\mathcal{F}_{t}] = M_{t}(x) $ and $M_{t}(x)$ is a non-negative continuous-time martingale that converges almost surely and in any $L^p$ for $p>1$ towards $\tilde{M}_{\infty}(x)$. This completes the proof of the proposition.%,  \begin{eqnarray} \label{convma} M_{t}(x) & \xrightarrow[t\to\infty]{a.s. \ L^p} & \tilde{M}_{\infty}(x).  \end{eqnarray}
\endproof

The process $x \in [0,1] \mapsto \tilde{M}_{\infty}(x)$ was used in \cite{BNS} to construct a fixed point to a recursive equation in distribution. In particular it is proved that $ x \mapsto \tilde{M}_{\infty}(x)$ is almost surely continuous.

\section{A geometric estimate}
In this section we establish  a rough control on the area of the largest and the smallest rectangle of $ \mathrm{Quad}(t)$. The reader can skip this part on first reading. For $t>0$, let $I_t := \inf \lambda(R)$ and $S_t := \sup \lambda(R)$ where the infimum and supremum run over all the rectangles $R \in \mathrm{Quad}(t)$. We will roughly prove that $t^{-4+o(1)} \leq I_{t}  $ and $S_{t} \leq t^{-1+o(1)}$. The formal statement is the following:
\begin{lemma} \label{area} For every $ \varepsilon >0$  we have

 \begin{eqnarray*}-4- \varepsilon < \liminf_{t\to \infty} \frac{\log\big(I_t)}{\log(t)} \ \ \mbox{a.s.}\qquad \mbox{ and } \qquad  \limsup_{t \to \infty}\frac{\log\big(P( S_{t} > t^{-1+\varepsilon})\big)}{\log(t)} = -\infty.\end{eqnarray*}
\end{lemma}
\proof \textsc{Lower bound.} Let $(x_i)_{1 \leq i \leq n}$ and $(y_i)_{1 \leq i \leq n}$ be the coordinates of the points of $ \Pi$ that occur before time $t$. By standard properties of Poisson point processes, conditionally on $n$, $(x_{i})$ and $(y_{i})$ are independent sequences of $n$ i.i.d.\,uniform variables over $[0,1]$.  A simple geometric argument (see Fig.\,\ref{min} below) shows that   \begin{eqnarray*}
I_t &\geq& \min_{\begin{subarray}{c} i \ne j\\1 \leq i,j \leq n \end{subarray} }|x_i-x_j| \cdot \min_{\begin{subarray}{c} i \ne j\\1 \leq i,j \leq n \end{subarray}}|y_i-y_j|.   \end{eqnarray*} 
\begin{figure}[!h]
 \begin{center}
 \includegraphics[width=6cm]{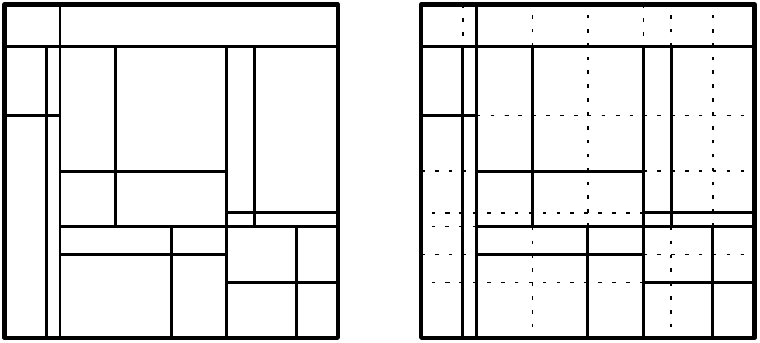}
 \caption{Illustration of the lower bound. \label{min}}
 \end{center}
 \end{figure}

Let $ \varepsilon >0$. By classical results on the uniform sieve of the interval $[0,1]$, if $x_{1}, x_{2}, ... $ are i.i.d.\,uniform points over $[0,1]$ then $\min \{|x_{i}-x_{j}| : 1 \leq i,j \leq n, i \ne j\}$ is asymptotically larger than $n^{-2- \varepsilon}$ a.s.\,. Indeed we have 
 \begin{eqnarray*}P\left( x_{n} \in \bigcup_{i=1}^{n-1} [x_{i}-n^{-2- \varepsilon}, x_{i}+ n^{-2- \varepsilon}] \right) &\leq& 2n^{-1- \varepsilon},  \end{eqnarray*} and an application of Borel Cantelli's lemma proves the claim. Since eventually $t/2 \leq n \leq 2t$ and $t \mapsto I_t$ is decreasing we almost surely have $I_t \geq t^{-4-2 \varepsilon}$ eventually.

 \textsc{Upper bound.} We use a common technique in fragmentation theory: the tagged particle. Assume that independently of the quadtree process $( \mathrm{Quad}(t))_{t \geq 0}$ we are given an independent variable $(U,V)$ uniformly distributed over $[0,1]^2$. The rectangle $ R^\bullet _{t}\in \mathrm{Quad}(t)$ containing $(U,V)$ is called the ``tagged rectangle'' at time $t$. The distribution of $(R^\bullet_{t})_{t \geq 0}$ is  equivalently described as follows. We start with $ R^\bullet_{0} := [0,1]^2$ and define the process $R_{t}^\bullet$ iteratively: when a splitting occurs at time $\tau$ inside the tagged rectangle $R^\bullet_{\tau^-}$, then $R_{\tau}^\bullet$ is one of the four subrectangles of $R_{\tau^-}^\bullet$ chosen proportionally to its two-dimensional Lebesgue measure. 
 
 It is clear from the above construction that the tagged rectangle at generation $n$ has a two-dimensional Lebesgue measure  which is distributed according to 
 $$ \prod_{i=1}^{2n} \overline{U}_{i},$$
 where $ \overline{U}_{i}$  are independent identically distributed variables with density $2m \mathbf{1}_{0<m<1}\mathrm{d}m$ (size-biased uniform over $[0,1]$). In particular we have $P(\prod_{i=1}^{2n}\overline{U}_{i} \geq \alpha ^n) \leq E[\overline{U}_{i}]^{2n}\alpha^{-n}=(4/9\alpha)^n$.   We now  turn to the study of the mass of the tagged rectangle at time $t>0$. For every $ \varepsilon >0$ we have 
   \begin{eqnarray*} P( \lambda( R^\bullet_{t}) > t^{-1+ 2\varepsilon}) & \leq & P( \mathrm{Gen}(R^\bullet_{t}) \leq t^{ \varepsilon}\,,\ \lambda( R^\bullet_{t}) > t^{-1 + 2 \varepsilon}) + P(\mathrm{Gen}(R^\bullet_{t}) > t^{ \varepsilon}\,,\ \lambda( R^\bullet_{t}) > t^{-1 + 2 \varepsilon}). \end{eqnarray*}
By our preceding remark, for large $t>0$, the second term of the right-hand side is bounded above  by $2^{- \lfloor t^\varepsilon \rfloor}$. For the first term, remark that if $ \lambda(R_{t}^\bullet) > t^{-1+2 \varepsilon}$ then for every $0 \leq s \leq t$, the intensity at which a particle falls inside $R_{s}^\bullet$ is larger than $t^{-1+2 \varepsilon}$, thus  by standard properties of exponential variables we have  $P( \mathrm{Gen}(R_{t}^\bullet)\leq t^{\varepsilon}\,,\ \lambda( R^\bullet_{t}) > t^{-1 + 2 \varepsilon})  \leq  P( \mathcal{P}(t^{2 \varepsilon}) \leq t^{ \varepsilon})$ where $ \mathcal{P}(t^{2 \varepsilon})$ is a Poisson distribution of mean $t ^{2 \varepsilon}$. Let us make this more precise. For $n \geq 0$, denote $\tilde{R}^\bullet_{n}$ the tagged rectangle at generation $n$. The rectangle $\tilde{R}_{n}^\bullet$ thus lives for an exponential time of parameter $\lambda ( \tilde{R}_{n}^\bullet)$ before it splits. We deduce that if $ \mathcal{E}_{0}, \mathcal{E}_{1},...$ is an i.i.d.\,sequence of exponential variables of parameter one which is also independent of $\lambda( \tilde{R}_{0}^\bullet), \lambda( \tilde{R}_{1}^\bullet), ... $ then  \begin{eqnarray*} P(\mathrm{Gen}(R_{t}^\bullet)\leq t^{\varepsilon}\,,\ \lambda( R^\bullet_{t}) > t^{-1 + 2 \varepsilon}) & = &P\left( \sum_{i=0}^{ \mathrm{Gen}(R_{t}^\bullet)} \lambda(\tilde{R}_{i}^\bullet)^{-1} \cdot  \mathcal{E}_{i} > t\,,\ \lambda( \tilde{R}^\bullet_{t}) > t^{-1+2 \varepsilon}\,,\ \mathrm{Gen}(R_{t}^\bullet) \leq t^ \varepsilon\right) \\ &\leq & 
P\left(  t^{1-2 \varepsilon} \sum_{i=0}^{ \mathrm{Gen}(R_{t}^\bullet)} \mathcal{E}_{i} > t \,,\ \mathrm{Gen}(R_{t}^\bullet) \leq t^{ \varepsilon}\right) \\ & \leq& P\left( \sum_{i=0}^{  \lfloor t^ \varepsilon \rfloor} \mathcal{E}_{i} > t^{ 2 \varepsilon} \right)\\ & = &
P\big( \mathcal{P}(t^{2 \varepsilon})\leq t^ \varepsilon\big).\end{eqnarray*}
The last probability being bounded above by $c^{-1}\exp({-ct^{ \varepsilon}})$ for some $c>0$. Gathering-up the pieces, there exists $c'>0$ such that we have $P( \lambda( R^\bullet_{t}) > t^{-1+ 2\varepsilon}) \leq c'^{-1} \exp(-c't^ \varepsilon)$ for every $t>0$. We then  use the tagged fragment to bound  $S_{t}$ from above. Notice that at any time $t>0$ the tagged fragment $R_{t}^\bullet$ has a probability $S_{t}$ of being the largest fragment, thus
  \begin{eqnarray*} P(\lambda(R_{t}^\bullet) > t^{-1+ 2 \varepsilon}) & \geq & P(S_{t}> t^{-1 + 2 \varepsilon}) t^{-1 + 2 \varepsilon},  \end{eqnarray*}which together with the previous bound easily completes the proof of the lemma.

%Using a similar argument as in the proof of the lower bound one can show that if we throw $\lfloor t^6 \rfloor$ independent and uniform points over $[0,1]^2$ then the probability that there exists one rectangle of $ \mathrm{Quad}(t)$ that does not contain one of these points is bounded above by $c't^{-2}$ for some $c'>0$. Since the law of a rectangle of $ \mathrm{Quad}(t)$ designated by a uniform point is that of the tagged rectangle we have for large $t>0$
 % \begin{eqnarray*} P( S_{t} > t^{-1 + 2 \varepsilon}) &\leq& c't^{-2} + t^{6} P( \lambda(R_{t}^\bullet) > t^{-1 + 2 \varepsilon}) \leq c' t^{-2}+ t^6c^{-1}e^{-ct^{ \varepsilon}}.  \end{eqnarray*}
%Since $E[S_{t}] \leq t^{-1+2 \varepsilon} + P(S_{t} > t^{-1+2 \varepsilon})$ the upper bound of the lemma easily follows.
%  Thus, using the fact that $t \mapsto S_{t}$ is decreasing, an application of Borel Cantelli's lemma shows that almost surely $S_{t} \leq t^{-1+2 \varepsilon}$ eventually, which completes the proof of the upper bound.
\endproof

\section{Proof of the main results}

\subsection{Proof of Theorem \ref{main}} \proof Let us first describe the main idea of the proof, which is similar to  \cite[Theorem 5]{BG04} and roughly speaking reduces to apply a law of large number after conditioning at a large time $t >0$. Fix $x \in[0,1]$ and let $T$ be much larger than $t$. Conditionally on $ \mathcal{F}_t$ the variable $N_T(x)$ is the sum of $N_t(x)$ independent contributions corresponding to the offsprings of the rectangles above $x$ at time $t$. By standard properties of the quadtree construction, the number of descendants of the rectangle $Q_t^i(x)$ inside $ \mathcal{Q}_T^i(x)$ is distributed as $N'_{ \lambda(Q_t^i(x))(T-t)}(x_{t}^i)$ where $N'_.(.)$ is an independent copy of the process $N_.(.)$.  Thus if for $x \in [0,1]$ and $t \geq 0$ we set $f(t,x) := E[N_t(x)]$, we have 

 \begin{eqnarray*} E\big[ N_T(x) \mid \mathcal{F}_t \big] &=&  \sum_{i \geq 1} f\big(\lambda(Q_t^i(x))(T-t),x_t^i\big)
. 
 %\label{variance}\mathrm{Var}\big(N_T(x) \mid \mathcal{F}_t\big) & = & \sum_{i \geq 1} \mathrm{Var}\big(N_{ \lambda(Q_t^i)(T-t)}\big).
 \end{eqnarray*} We now let $T \to \infty$ in the last display. Since for each rectangle $Q_{t}^i(x)$ of $ \mathcal{Q}_{t}(x)$ we have $(T-t)\lambda(Q_t^i(x)) \to \infty$ then $T^{-\beta}f(\lambda(Q_t^i(x))(T-t),x_t^i)$ tends to $K_{0} \lambda(Q_{t}^i)^\beta h(x_t^i)$ as $T \to \infty$. Henceforth we have
% By the choice of $t$ versus $T$ and the lower bound of Lemma \ref{area} we get that almost surely, for each rectangle $Q_{t}^i(x)$ of $ \mathcal{Q}_{t}(x)$ we have $(T-t)\lambda(Q_t^i(x)) \to \infty$ uniformly as $T \to \infty$ almost surely. By \eqref{CJ} we have $t^{-\beta}f(t,x) \to K_0 \cdot h(x)$ as $t \to \infty$.  Henceforth we have 
  \begin{eqnarray}&&\left|T^{-\beta}E\big[ N_T(x) \mid \mathcal{F}_t \big] - K_0\cdot \sum_{i \geq 1} \lambda(Q_t^i(x))^\beta  h(x_t^i) \right| \nonumber \\ &=& \left|T^{-\beta}E\big[ N_T(x) \mid \mathcal{F}_t \big]- K_0 \cdot M_t(x) \right| \quad\xrightarrow[T\to\infty]{} \quad0.  \label{lfgn} \end{eqnarray}
  The strategy of the proof is now clear: conditionally on $ \mathcal{F}_{t}$, by the law of large numbers, $T^{-\beta}N_{T}(x)$ will be very close to its (conditional) mean which is close to $K_{0} \cdot M_{t}(x)$ which converges towards $K_0\cdot \tilde{M}_\infty$. This will imply the theorem.

  To make this statement precise, and in particular get an almost sure convergence (a convergence in probability would be much easier to prove), we shall need the estimates on the expectation and the variance of the process $N_{t}(x)$ proved by   Broutin, Neininger and Sulzbach.  It follows from Proposition 11 and Theorem 5 in \cite{BNS} that there exist two constants $C>0$  and $\delta>0$ such that for every $t>0$ we have 
  \begin{eqnarray} \label{boundexp} \sup_{0 \leq x \leq 1} \big | t^{-\beta}E[N_t(x)]-K_0\cdot h(x) \big| & \leq & C t^{-\delta} \\
  \label{boundvar}  \sup_{0 \leq x \leq 1} \big(\mathrm{Var}(N_t(x)) \big) &\leq & C (t^{2\beta}+t), \end{eqnarray} the term $t$ appearing in the last line since the variance of $ N_{t}(x)$ is of order $t$ near $t=0$. We first make \eqref{lfgn} quantitative in $T$. Fix $\alpha \geq 5$ such that $\delta(\alpha-5) >1$ and for $t>0$, set $T:= t^\alpha$.  By the choice of $t$ versus $T$ and the lower bound of Lemma \ref{area} we get that almost surely, there exists a random time $\tau$ such that for $t \geq \tau$ we have $(T-t)\inf\{\lambda(Q_t^i(x)) : i \geq 1\} \geq t^{\alpha-5}$. Henceforth using the bound \eqref{boundexp}, we have for $t \geq \tau$
   \begin{eqnarray*} \left|(T-t)^{-\beta}E\big[ N_T(x) \mid \mathcal{F}_t \big]- K_0 \cdot M_t(x) \right| & \leq &  \sum_{i \geq 1} \Big| (T-t)^{-\beta}f\big( \lambda(Q_t^i(x))(T-t),x_t^i) - K_0\cdot \lambda(Q_t^i(x))^\beta h(x_t^i)\Big|  \\ & \leq &  C N_t(x) t^{-\delta(\alpha-5)}. \end{eqnarray*}
Since $N_t(x)$ is clearly less that the number of points fallen so far, the definition of $\alpha$ implies that   $N_t(x) t^{-\delta(\alpha-5)}$ goes to $0$ almost surely. Since $M_t(x)$ is  almost surely bounded,  we proved that with our choice of $\alpha$ we have   $|T^{-\beta}E[N_T(x) \mid \mathcal{F}_t]- K_0 \cdot M_t(x)| \to 0$ almost surely as $t \to \infty$ and using Proposition \ref{propmartingale} if follows that  \begin{eqnarray} \label{lfgn2}  \big|T^{-\beta}E[N_T(x) \mid \mathcal{F}_t]- K_0 \cdot \tilde{M}_{\infty}(x)\big| & \xrightarrow[t\to\infty]{a.s.}& 0.  \end{eqnarray}

  Recall that conditionally on $ \mathcal{F}_{t}$ the contributions of each rectangle $Q_{t}^i(x)$ to $N_{T}(x)$ are independent, thus we have 
   \begin{eqnarray*}T^{-2\beta}E \left[ \big( N_{T}(x)- E[N_{T}(x) \mid F_{t}]\big)^2 \mid \mathcal{F}_{t}\right] &=&T^{-2\beta}\sum_{i \geq 1} \mathrm{Var}\big(N_{ \lambda(Q_t^i(x))(T-t)}(x_t^i)\big) \\ & \leq &  C \Big( \sum_{i \geq 1} \lambda(Q_t^i(x))^{2\beta}+ T^{1-2\beta} \sum_{i \geq 1} \lambda(Q_{t}^i(x))\Big)\\
& \leq& C \big(  S_t^{2\beta-1} + T^{1-2\beta}\big) \sum_{R \in \mathrm{Quad}(t)} \lambda(R)\\ &=&   C(  S_t^{2\beta-1} + T^{1-2\beta}) , \end{eqnarray*}
where we used the bound \eqref{boundvar} to reach the second line and the fact that $2\beta >1$ to go from the third to the last line. Let $\varepsilon >0$. Applying the standard Markov inequality conditionally on $ \mathcal{F}_{t}$ we obtain
  \begin{eqnarray} \label{equa} P \Big( \left.T^{-\beta}\big|N_T(x)-E[N_T(x) \mid \mathcal{F}_t]\big | \geq \varepsilon \ \right|\ \mathcal{F}_t\Big) & \leq & C \varepsilon^{-2}\big(S_t^{2\beta-1} + T^{1-2\beta}\big).  \end{eqnarray}
  We now take $T_{k}:=(1+\eta)^k$ and $t_{k}:=(1+\eta)^{k/\alpha}$ for $k=1,2,3...$ and $\eta >0$. Since $E[S_{t_{k}}^{2\beta-1}]$ is less than $t_{k}^{(2\beta-1)(\varepsilon-1)} + P(S_{t_{k}}> t_{k}^{-1+ \varepsilon})$, using Lemma \ref{area} and \eqref{equa} we see that $P(T_{k}^{-\beta}|N_{T_{k}}(x)-E[N_{T_{k}}(x) \mid \mathcal{F}_{t_{k}}] | \geq \varepsilon)$ is summable in $k$.  Applying Borel-Cantelli's lemma and using  \eqref{lfgn2} we deduce that for every $\eta>0$ we have the following almost sure convergence  \begin{eqnarray*}(1+\eta)^{-k\beta} N_{(1+\eta)^k}(x) & \xrightarrow[k\to\infty]{a.s.} &K_{0}\cdot \tilde{M}_{\infty}(x).  \end{eqnarray*}To extend this result to the whole process we use the fact that $t \mapsto N_{t}(x)$ is increasing in $t$ which implies
   \begin{eqnarray*} (1+\eta)^{-(k+1)\beta} N_{(1+\eta)^k}(x) \leq s^{-\beta} N_{s}(x) \leq (1+\eta)^{-k\beta} N_{(1+\eta)^{k+1}}(x)  \end{eqnarray*} for every $(1+\eta)^k \leq s \leq (1+\eta)^{k+1}$ and $k \geq 1$. Since this holds for any $\eta>0$ we easily deduce that $t^{-\beta}N_{t}(x)$ almost surely converges towards $K_{0} \cdot \tilde{M}_{\infty}(x)$. This completes the proof of Theorem \ref{main}.\endproof 
  
 \subsection{Proof of Corollary \ref{coro}}
 \proof(Sketch) Theorem \ref{main} implies the convergence of the finite dimensional marginals of
 $t^{-\beta}N_{t}(.)$ towards those of $ K_{0}\cdot \tilde{M}_{\infty}(.)$ in probability: For any $0 \leq x_{1}, ... , x_{k} \leq 1$ we have 
  \begin{eqnarray} \label{fidi} t^{-\beta}\big( N_{t}(x_{i}) \big)_{1 \leq i \leq k} & \xrightarrow[t\to\infty]{(P)}& K_{0} \big( \tilde{M}_{\infty}(x_{i})\big)_{1\leq i \leq k}.  \end{eqnarray} Furthermore Theorem 1 of \cite{BNS} provides the tightness of the processes $t^{-\beta}(N_{t}(.))$ for the uniform metric:  for every $ \varepsilon>0$ there exists $\eta >0$ such that for $t>0$ large enough we have 
    \begin{eqnarray} \label{tight}P(\omega_{ t^{-\beta}N_{t}(.)}( \eta)  \leq  \varepsilon) \geq 1- \varepsilon,  \end{eqnarray} where $\omega_{g}(\eta) = \sup\{ |g(x)-g(y)|\,,\ |x-y| \leq \eta\}$ is the modulus of continuity of the function $g$. Recalling that $ x \mapsto \tilde{M}_{\infty}(x)$ is almost surely continuous, we can combine \eqref{fidi} and \eqref{tight} to get that $t^{-\beta}N_{t}(.)$ converges in probability for the $L^\infty$ metric towards $ \tilde{M}_{\infty}(.)$.  We leave the details to the reader.\endproof 
    
\noindent \textbf{Open Question.} It is believable that the convergence of Corollary \ref{coro} actually holds almost surely, that is 
 \begin{eqnarray*}\big( t^{-\beta} N_{t}(x)\big)_{x \in [0,1]}& \xrightarrow[t\to\infty]{a.s.} & K_{0} \cdot \big(\tilde{M}_{\infty}(x)\big)_{x \in [0,1]},  \end{eqnarray*} for the uniform metric $\|.\|_\infty$.

 \section{Fragmentation process with parameter} \label{comments}
 
 In this section we comment at an informal level on the strategy adopted in this work and on possible extensions of our techniques.
 
%The overall strategy adopted in this work bears many similarities with the theory of self-similar fragmentation theory. Indeed, in \cite{CJ10} it was shown that the masses of the rectangles $\{Q_{t}^i(x)\}$ is \emph{not} a fragmentation process, although some information was derived using fragmentation theory. More precisely it was shown that the number of rectangles above a uniform point $U$ at time $t>0$ has the same expectation as the number of fragments in a related fragmentation process where, roughly speaking the positions are forgotten.\\

\paragraph{Fragmentation theory.} Let us  briefly recall some basics about fragmentation theory. We stick to a very simple case for sake of clarity.  For more details, we refer to \cite{Ber06}. To define a  self-similar fragmentation process\footnote{binary, without erosion and with dislocation measure of mass one} $ \mathscr{F}$ we need one input:  a probability measure $\nu$  on $\{ (s_1, s_2) : s_{1}\geq s_{2} >0 \textrm{ and } s_{1}+s_{2} \leq 1\}$. 
 The process $ \mathscr{F}$ with values in the set $\mathcal{S}^\downarrow = \{ (s_1, s_2, \dots) : s_{1} \geq s_2 \geq \dots  \geq 0 \textrm{ and } \sum_{i} s_{i} \leq 1\}$ is informally characterized as follows: if at time $t$ we have $\mathscr{F}(t)=(s_{1}(t), s_{2}(t), \dots )$, then for every $i\geq 1$, the $i$-th ``particle'' of mass $s_{i}(t)$ lives an exponential time with parameter $s_{i}(t)$  before splitting into two particles of masses $r_{1}s_{i}(t)$ and $r_{2}s_{i}(t)$, where $(r_{1},r_{2})$ has been sampled from $\nu$ independently of the past and of the other particles. In other words, each particle undergoes a self-similar fragmentation with time rescaled by its mass. It is  classical that under mild assumption there exists a unique $\beta \in (0,1]$ (called the Malthusian exponent) such that 
  \begin{eqnarray} \label{malthus} \int \mathrm{d}\nu(s_{1},s_{2}) \ s_{1}^\beta + s_{2}^\beta &=&1,  \end{eqnarray} and that $ M_{t} := \sum s_{i}(t)^\beta$ is a continuous-time non-negative martingale which plays a central role in the asymptotic behavior of the fragmentation process, see \cite{Ber06,BG04}.

\paragraph{Parametrized fragmentation.} In the problem of the partial match query, one can think of the rectangles above the point $x$ at time $t>0$ as a \emph{fragmentation process where the particles have an additional parameter}, in our case the position $x_{t}^i \in [0,1]$. This leads us to extend the notion of dislocation measure and to define a fragmentation process with parameter: a parametrized (binary) dislocation probability is a collection $ \boldsymbol{\nu} = (\nu_{x})_{x \in [0,1]}$ such that for every $x \in [0,1]$, $\nu_{x}$ is a probability measure on $$\{(s_{1},x_2,s_{2},x_{2}) \in [0,1]^4 : s_{1} \geq s_{1} \ \mathrm{ and }\  s_{1} + s_{2} \leq 1\}.$$ A parametrized fragmentation process $ \mathcal{F}$ with dislocation measure $ \boldsymbol{\nu}$ is then a process with values in $\{ (s_1,x_{1}, s_2, x_{2}, ... \dots)  \in [0,1]^ \mathbb{N}: s_{1} \geq s_2 \geq \dots  \geq 0 \textrm{ and } \sum_{i} s_{i} \leq 1\}$ whose evolution is informally described as follows: We start with a particle of mass $1$ given with a position $ x \in [0,1]$. If  $ \mathcal{F}(t) = (s_{1}(t), x_{1}(t), s_{2}(t), x_{2}(t), \dots)$ then for every $i\geq 1$, the $i$-th ``particle'' of mass $s_{i}(t)$ with position $x_{i}(t)$ lives an exponential time with parameter $s_{i}(t)$  before splitting into two particles of masses $r_{1}s_{i}(t)$ and $r_{2}s_{i}(t)$ with respective positions $y_{1}$ and $y_{2}$, where $(r_{1},y_{1},r_{2},y_{2})$ has been sampled from $\nu_{x_{i}(t)}$ independently of the past and of the other particles. In this setting, \eqref{malthus} is replaced by the following assumption: $(H)$ There exists $\beta \in [0,1]$ and $h : x \in [0,1] \mapsto h(x) \in \mathbb{R}_{+}$ such that for every $x \in [0,1]$ we have 
  \begin{eqnarray} \label{malthusbis} \int \mathrm{d} \nu_{x}(s_{1},x_{1},s_{2},x_{2}) \big(s_{1}^\beta h(x_{1})+ s_{2}^\beta h(x_{2})\big) &=& h(x).  \end{eqnarray}
  Then under this assumption the process $M_{t}(x) := \sum s_{i}(t)^\beta h(x_{i}(t))$ is a continuous-time non-negative martingale playing the role of the Malthusian martingale. It is believable that substantial parts of self-similar fragmentations theory can be adapted to this parametrized case. \medskip
  
\paragraph{Examples.} \textsc{Partial Match queries in Quadtree.} Within this formalism the process of the masses of the rectangles of $ \mathcal{Q}_t(x)$ is a parametrized fragmentation process starting with a single particle of mass $1$ and parameter $x$. Its parametrized dislocation measure $ \boldsymbol{\nu^{ \mathrm{quad}}}$ is given by 
 \begin{eqnarray*}
&&\int \mathrm{d}\nu^{ \mathrm{quad}}_x(s_1,x_1,s_2,x_2) f(s_1,x_1,s_2,x_2)\\  &=& \iint_{[0,1]^2}\mathrm{d}u\textrm{d} v  \Bigg( \textbf{1}_{x < u} f \left ( uv, \frac{x}{u},u(1-v),\frac{x}{u} \right ) +  \textbf{1}_{x > u} f \left ( (1-u)v, \frac{x-u}{1-u}, (1-u)(1-v),\frac{x-u}{1-u} \right ) \Bigg),  \end{eqnarray*}
 for every $x \in [0,1]$ and  every Borel function $f : [0,1]^4 \to \mathbb{R}_+$. In particular $\mathrm{d}\nu_{ \mathrm{quad}}(s_1,x_1,s_2,x_2)$-almost surely we have $x_1=x_2$ and Hypothesis $(H)$ is fulfilled with $\beta= \frac{ \sqrt{17}-3}{2}$ and $h(x)=(x(1-x))^{\beta/2}$. \bigskip

  To conclude, besides the application of the method to the problem of partial match queries in higher-dimensional random quadtrees or in random $k$-d trees, we present another setup taken from \cite{CLG10} where the concept of  ``parametrized fragmentation'' could be applied (although the results there only rely on the standard fragmentation theory). \bigskip

  \textsc{Random chords.} We  recall the random chord construction of \cite{CLG10}. We consider a sequence $U_1,V_1,U_2,V_2,\ldots$ of independent
random variables, which are uniformly distributed over the unit circle $\sun$. We then construct inductively
a sequence $L_1,L_2,\ldots$ of random closed subsets of the (closed) unit disk $\overline{\mathbb D}$.
To begin with, $L_1$ just consists of the chord with endpoints $U_1$, and $V_1$, which we denote 
by $[U_1V_1]$. Then at step $n+1$, we consider two cases. Either the chord $[U_{n+1}V_{n+1}]$
intersects $L_n$, and we put $L_{n+1}:=L_n$. Or the chord $[U_{n+1}V_{n+1}]$ does not
intersect $L_n$, and we put $L_{n+1}:=L_n \cup [U_{n+1}V_{n+1}]$. Thus, for every integer $n\geq 1$, 
$L_n$ is a disjoint union of random chords.  If $x,y \in \mathbb{S}_{1}$ then one defines the fragments separating $x$ from $y$ as the connected components of $ \overline{\mathbb{D}} \backslash L_{n}$ intersecting $[x,y]$.

\begin{figure}[!h]
 \begin{center}
 \includegraphics[width=13cm]{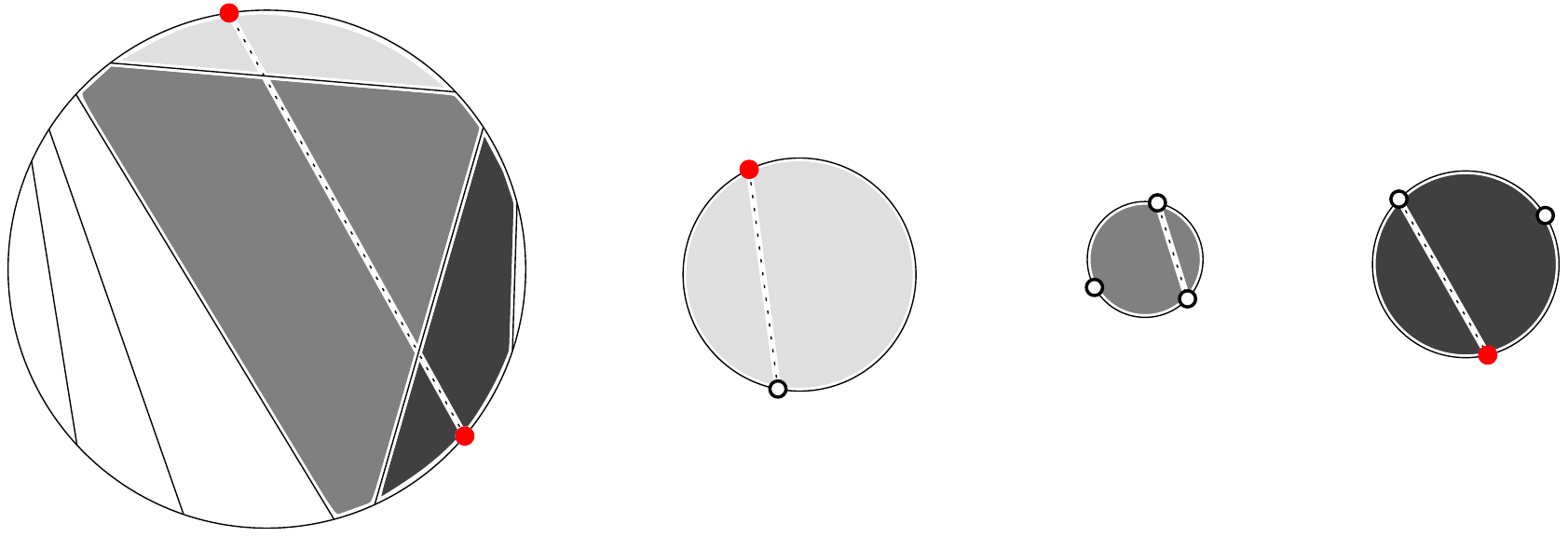}
 \caption{The fragments separating two points. \label{separating}}
 \end{center}
 \end{figure} 
 
 After contracting the chords of $L_{n}$, each fragment $F$ separating $x$ from $y$ can be seen as a particle with two distinguished points (see Fig. \ref{separating}) whose mass corresponds to the one-dimensional Lebesgue measure of $F \cap \mathbb{S}_{1}$, see \cite{CLG10}. The ``position'' or parameter of each particle is then the relative positions of its two distinguished points in $[0,1]$.  It was shown in \cite{CLG10} that if $U,V$ are independent and uniformly distributed over $ \mathbb{S}_{1}$ then the fragments separating $U$ from $V$ form (in a proper continuous time parametrization) a fragmentation process. In the case when $x,y \in \mathbb{S}_{1}$ are fixed, the process of fragments separating $x$ from $y$ (in a proper continuous time parameterization) can be seen as a parametrized fragmentation process with assumption $(H)$ fulfilled with $$\beta = \frac{\sqrt{17}-3}{2} \quad  \mbox{and} \quad  h(x) = (x(1-x))^\beta,$$ which is equivalent to equation $(14)$ of \cite{CLG10}. % The convergence of Proposition implies the convergence in probability of the process $t^{-\beta}N_{t}(.)$ towards $ K_{0} \cdot\tilde{M}_{\infty}(.)$ in the sense of finite dimensional marginals. Furthermore it follows from \cite[Theorem 1]{BNS} that the sequence of laws of the processes $(t^{-\beta}N_{t}(.))_{t \geq 0}$  \endproof

\bibliographystyle{abbrv}

\noindent
D\'epartement de Math\'ematiques et Applications 
\\
\'Ecole Normale Sup\'erieure, 45 rue d'Ulm
\\
75230 Paris cedex 05, France

\medbreak

\noindent
nicolas.curien@ens.fr

 \end{document}